\documentclass[11pt]{article}
\usepackage[american]{babel}
\usepackage[utf8]{inputenc}
\usepackage[T1]{fontenc}
\usepackage{lmodern}
\usepackage[dvips]{graphicx}

\usepackage{amsfonts}
\usepackage{amsmath}
\usepackage{enumerate}
\begin{document}

\title{Some notes on the impact of Lagrange's memoir  
``On the construction of geographical maps"}
\author{Athanase Papadopoulos\thanks{Institut de Recherche Mathématique Avancée
(Université de Strasbourg et CNRS)
7 rue René Descartes
67084 Strasbourg Cedex France
email : athanase.papadopoulos@math.unistra.fr}}

\date{}
\maketitle

      \begin{abstract} 
 These are notes on the impact of Lagrange's memoir on the construction of geographical maps. We mention the relations of some ideas and questions introduced in this memoir with other notions that appeared later in the works of several mathematicians, including in particular Chebyshev (19th c.) and Darboux (19th-20th c.), two mathematicians who were particularly interested in geography. The final version of this paper appears in the book: Mathematical Geography in the Eighteenth Century: Euler, Lagrange and Lambert (ed. R. Caddeo and A. Papadopoulos), Springer, 2022.
     \end{abstract}
     
     \bigskip
\bigskip

          \section{Introduction} In this chapter, I shall mention the impact of some ideas and questions introduced by Lagrange in his memoir \emph{Sur la construction des cartes g\'eographiques} (On the construction of geographical maps) \cite{Lagrange-Construction} on the works of some later authors including Darboux\index{Darboux, Gaston} (\S  \ref{s:Darboux}),   Chebyshev (\S  \ref{s:Chebyshev}) and others (\S  \ref{s:Chebyshev} and  \ref{s:later}). The last section (\S \ref{s:modern}) contains some notes on the use of Lagrange's ideas in the actual construction of geographical maps.

          \section{Lagrange's memoir in  the work of Darboux}\label{s:Darboux}
 
Lagrange,  in his memoir \emph{Sur la construction des cartes g\'eographiques}, motivated by works of Lambert and Euler, addressed the following problem: 

\emph{To find all the conformal projections from the sphere (and more generally from a surface of revolution) onto the Euclidean plane by which the meridians and parallels are sent to circles.}  

This problem of Lagrange is mentioned several times by Darboux,\index{Darboux, Gaston} in his
\emph{Leçons sur la th\'eorie g\'en\'erale des surfaces et les applications g\'eom\'etriques du calcul infinit\'esimal} \cite{Darboux-Lecons}, and in particular in
Chapter IV of Book II of Part I, which concerns conformal maps between surfaces. In that chapter, Darboux,\index{Darboux, Gaston} after having studied in the preceding chapters the general methods of isothermal coordinates,\index{isothermal coordinates} applies them to the resolution of  problems posed by the construction of geographical maps. He  considers Lagrange's problem (referring to the latter as the first mathematician who addressed it)  in the setting of the spheroid. Let us recall that in this setting, a \emph{meridian} is an ellipse (or half-ellipse) whose rotation about an axis generates the spheroid, and a \emph{parallel} is the orbit of a point on this ellipse by this rotation. The meridians are geodesics of the resulting surface but the parallels are not. (This holds in the general case of a surface of revolution in Euclidean  space). Darboux\index{Darboux, Gaston} makes the remark that if one of the two families (the meridians and the parallels) is represented by circles, then the same holds for the second family. Thus, Lagrange's conditions (that both families are represented by circles) are redundant. He discusses at length the various solutions to this problem using geometric transformations, see  \cite[p. 236-241]{Darboux-Lecons}.
In the same chapter, Darboux reviews a method of Gauss on the conformal representation on a sphere of a region of the Earth whose figure is assumed to be spheroidal  such that the similarity ratio is constant along a meridian. Applying this method to the map of France, Darboux\index{Darboux, Gaston} finds that the variation of the similarity ratio, on the whole extent of the map, is only of the order of $\frac{1}{400.000}$ of its actual value.
 
  Darboux gives a modern formulation and a detailed modern proof of the main problem stated by Lagrange in his memoir and which we recalled at the beginning of this section. In Darboux' formulation,\index{Darboux, Gaston} the solution consists of the following steps:

\begin{enumerate}
\item  Perform a stereographic projection on the equator plane, with the North pole as center. This gives a first map, call it A.

\item Transform the map A by a mapping defined by the equations
\[\rho^c=\rho', \hskip.2in c\omega=\omega'
\]
where $(\rho,\omega)$ and $(\rho',\omega')$ are polar coordinates and where the pole is taken to be center of the sphere
(Darboux calls this a Lambert transformation).\index{Lambert transformation} The constant $c$ is what Lagrange calls the \emph{exponent} of the projection.\index{exponent of projection}
 This gives a second map, call it B.
\item Apply an arbitrary inversion, with pole in the plane of the equator. This gives the needed map.
\end{enumerate}

Let me mention now another problem of Lagrange to which he was led in the same memoir, and of which  Darboux\index{Darboux, Gaston} gave a solution. This is a problem in plane geometry, whose statement is the following:

 \emph{Given three points $R$, $R'$, $R''$, to construct on a fixed  basis $AB$ three triangles whose vertices are the given points and such that:
\begin{enumerate}
\item the differences of the angles at the vertices $BRA$, $BR'A$, $BR''A$ are given; 
\item the ratios of the sides containing these angles, that is, $\frac{RB}{RA}$, $\frac{R'B}{R'A}$, $\frac{R''B}{R''A}$, are given.
\end{enumerate}
}

Lagrange says that this problem seems to him quite difficult to solve by geometry, and that he did not try the algebraic solution because this seemed to him useless, unless one can reduce it to an easy construction. He then sketches a solution of the problem using the analytical methods he developed in the first part of the memoir.

Darboux, some 130 years after Lagrange, published an article titled \emph{Sur un probl\`eme pos\'e par Lagrange} (On a problem posed by Lagrange) \cite{Darboux-Lagrange} in which he gave a solution of this problem using the new geometrical methods that were available at his time, after transforming it into a problem of finding an inversion in the plane which transforms one triangle into another one. Darboux\index{Darboux, Gaston} writes that the advances made in geometrical methods allow  now to solve easily and elegantly the problem that Lagrange  indicated, and that this problem is reduced to the following: 

\emph{Given a triangle $ABC$ in the plane, to find an inversion that turns it into a triangle equal to another given triangle $A_0B_0C_0$.}  

Darboux then gives a concise geometric and complete solution to this problem. (The same solution by Darboux is contained in Part I, \S 133, Chapter IV of his  \emph{Leçons sur la th\'eorie g\'en\'erale des surfaces et les applications g\'eom\'etriques du calcul infinit\'esimal}, p. 241-243.)

\section{Chebyshev's theorem and developments}
\label{s:Chebyshev} 

At the end of \S 2 of his memoir \cite{Lagrange-Construction}, Lagrange gives a formula for the infinitesimal dilatation\index{infinitesimal dilatation} at a point of an angle-preserving map between (a subset of) the sphere and the plane. The notation is the following:
$s,t$ are the curvilinear coordinates of a point on the surface of the Earth (assimilated to a sphere of radius 1), and 
$x,y$ are the rectangular coordinates of its image in the plane. The differential $ds$ denotes the difference of the two arcs of meridian passing through the infinitesimally close points on the sphere and $qdt$ denotes the arc of parallel contained by these two meridians. 
Then, $(s+ds, t+dt)$ and $(x+dx, y+dy)$ are respectively the curvilinear coordinates of a point on the sphere which is infinitesimally close to $(s,t)$ and those of its image by the map. With this, Lagrange writes the following formula for the \emph{infinitesimal dilatation}\index{infinitesimal dilatation}  at the given point on the sphere:
$$
m=\sqrt{ \frac{ dx^2+dy^2}{ ds^2 + q^2 dt^2}}.
$$

 This formula has been used by  P. L. Chebyshev\index{Chebyshev, Pafnuty Lvovich} as a starting point for his investigations on the construction of geographical maps that minimize distortion. More precisely, using the above notation, in two papers, carrying the same titles as those of Lagrange, Chebyshev addressed the question of determining the maps that minimize the deviation of the infinitesimal dilatation\index{infinitesimal dilatation}  $m$ from its integral over the region considered on the sphere. He used other formulae established by Lagrange in his paper and he obtained a new result, namely, that among all the conformal representations of a subset of the sphere onto the plane, the representation that minimizes distortion is the one for which this distortion is constant on the boundary of the region,  see the memoirs \cite{Chebyshev1} and \cite{Chebyshev2}. Chebyshev\index{Chebyshev, Pafnuty Lvovich} also made a relation between this study and the Laplace equation.
   
 Since in the previous section we talked about Darboux'\index{Darboux, Gaston} work on geography, let us mention here that  
in 1911, Darboux published a memoir, whose title is the same as those of Lagrange and Chebyshev (On the construction 
of geographical maps), in which he presented a proof of Chebyshev's
theorem, see \cite{Darboux-construction}. As a matter of fact, in his paper, Darboux attributes to
Chebyshev\index{Chebyshev, Pafnuty Lvovich} a more general result, valid for any surface of
revolution.

Let us also mention the work of D. A. Grav\'e,\index{Grav\'e, Dimitri Aleksandrovich} a student of Chebyshev. 
Whereas Lagrange, in his paper \cite{Lagrange-Construction}, considered the problem of finding all the angle-preserving maps for which the meridians and the parallels are sent to circles or straight lines,  Grav\'e worked on the problem of characterizing the area-preserving maps for which the meridians and the parallels are sent to circles or straight lines. He gave a complete solution of this problem in a paper which is also titled \emph{Sur la construction des cartes g\'eographiques}, see \cite{Grave1869}. Furthermore,  in 1894,  Grav\'e\index{Grav\'e, Dimitri Aleksandrovich}  presented  an outline of a proof of Chebyshev's theorem (which the latter has only sketched), at the annual meeting of the \emph{Association fran\c caise pour l'avancement des sciences} which took place in Caen.\footnote{There is a written note by Grav\'e titled \emph{Sur une question de Tch\'ebychef} (On a question of Chebyshev)  in the 1895 publications of the Association.
Regarding the many possible transliteration of Chebyshev's name from the cyrillic into Latin characters,  Grav\'e, uses  \emph{Tchebichef}, the one that was used by Chebyshev himself when he signed the (numerous) papers he wrote in French.   In fact, this is the unique way in which Chebyshev wanted his name to be transliterated. M. d'Ocagne published a note on this matter in the \emph{Bulletin des sciences math\'ematiques}, titled ``On the spelling of the name of Tchebichef" \cite{Ocagne1931},  in which he writes: ``Because of the many variants of the spelling of the name Tchebichef which one can find in the contemporary mathematical publications, let me recall the following fact which I have already pointed out.
During his last stay in Paris, in May 1893, the illustrious Russian geometer accepted to entrust me with the task of writing a complete description (which had never been given even in Russian) of his intriguing arithmetic machine, a description which he reviewed himself with the  greatest care.
On this occasion, he personally prescribed the French spelling of his name under the form `Tchebichef', formulating the wish that it is exclusively adopted in all his writings written with our characters. It is because he recommended me to ensure his \emph{desideratum} regarding this point, from which common usage deviates,  that the preceding lines have been written." Obviously, Chebyshev's  wish was not fulfilled.} In 1896, Grav\'e published in Russian a complete proof of Chebyshev's theorem, and in 1911 he published a paper in French  titled \emph{Sur un th\'eor\`eme de Tch\'ebychef g\'en\'eralis\'e} (on a generalized theorem of Chebyshev), in which he gave a detailed proof of a slightly more general result, valid  not only for the sphere, but for an arbitrary surface whose curvature does not change sign (see \cite{Grave1911}).

  Finally, we note that a modern  proof of Chebyshev's  theorem, which is based upon his ideas,  was given about a century after Chebyshev found it,  by John Milnor\index{Milnor, John} \cite{Milnor}, who also pointed out  further developments, highlighting the case where the region of the sphere which is mapped is geodesically convex.
   Milnor writes in his paper: ``This result has been available for more than a hundred years, but to my knowledge it has never been used by actual map makers."

\section{Lagrange's memoir in the modern literature}\label{s:later}

We start this section by recalling the idea of the Schwarzian derivative.\index{Schwarzian derivative} This is a function associated with a smooth function of a real variable, or a holomorphic function of a complex variable, which measures the deviation of such a function from being a M\"obius transformation. 
A more precise definition considers the Schwarzian derivative as a 1-cocycle on the group of diffeomorphisms of the projective space $\mathbb{RP}^1$ with coefficients in the space of quadratic differentials  (holomorphic quadratic differential forms) and kernel in $\mathrm{PSL}(2,\mathbb{R})$; the  Schwarzian derivative is the only such cocycle, see \cite{OT}. The clearest geometric definition is probably  that given by Thurston in \cite{Th}, in which a Schwarzian derivative becomes simply a quadratic differential. His aim in his paper is to study a topology on the space of conformal maps from the unit disc to simply connected domains of the complex plane; such a topology is induced by the topology of uniform convergence of Schwarzian derivatives. In the introduction of his paper, Thurston makes an analogy between the notion of Schwarzian derivative and that of curvature in differential geometry.  

V. Ovsienko and S. Tabachnikov notice in  \cite{O, OT} that the Schwarzian derivative was already introduced by Lagrange in his paper on the construction of geographical maps.  It is also reported in \cite{OT} that Hermann Schwarz, who introduced Schwarzian derivatives,\index{Schwarzian derivative} declared that this notion is present at least implicitly in Lagrange's memoir. According to an explanation given to the author of this chapter by Ovsienko, it is in \S 10 of Lagrange's memoir, starting with  the definition of $\varphi$ and $\Phi$, that this notion occurs. Indeed, Lagrange works with the quotients $\varphi''/\varphi$ and $\Phi''/\Phi$ and he compares them, and two functions, $f$ and $F$ appear in this context.  A small computation of  $\varphi''/\varphi$ in terms of $f$ and $\Phi''/\Phi$ in terms of $F$ shows that they are equal precisely to $1/2 S(f)$ and $1/2S(F)$, where $S$ denotes the Schwarzian derivative. The result is that Lagrange expresses the property of sending circles to circles in terms of having equal Schwarzian derivatives.

Another relation of Lagrange's memoir with modern works concerns the point of view of metric geometry. The problem of finding maps from the sphere to the plane in which circles are sent to straight lines, which is a particular case of the general problem Lagrange addressed in his memoir on geography, has been considered by V. S. Matveev and others as an early formulation of the general problem of finding maps between surfaces that send geodesics to geodesics, and of the theory of geodesically
equivalent metrics, see the papers \cite{Mat, BV}, 
 and the older paper by Beltrami \cite{Beltrami} in which the latter says that his work is motivated by the question of constructing geographical maps.

\section{On Lagrange's work and practical cartography}\label{s:modern}

It is not clear to what extent Lagrange's work has been used in practical cartography, but there are indications for the fact that some attempts were made for that.
In a letter to Lagrange dated February 14, 1782, Laplace writes  \cite[Vol. XIV, p. 111]{Lagrange-Oeuvres}:  ``[\ldots] Your two memoirs on the construction of geographical maps gave me as much pleasure. Above all, I admired the elegant manner with which you extracted from the general solution of the problem the case where the meridian[s] and the parallels are represented by circles. Besides, your analysis has the merit of being useful in the practice of constructing particular maps, and I have engaged one of my friends, who just announced a big atlas, to use it."

Around the beginning of the 20th century,  Jules de Schokalsky,\index{Schokalsky@de Schokalsky, Jules} a colonel of the Imperial Russian Navy and secretary of the Physics Section of the Imperial Russian Geographical Society, was entrusted with the publication of an atlas of universal geography, containing, among other things, a series of maps of European Russia on a scale of 1: 1,640,000.
From N. de Zinger's  article \emph{La projection de Lagrange appliqu\'ee \`a la carte de la Russie d'Europe} (The Lagrange projection applied to the map of European Europe)  \cite{Zinger}, we learn that it was a so-called ``Lagrange projection"  that was chosen, with some convenient constants. Even though it is not clear to what part of Lagrange's memoir the author refers, this indicates at least that Lagrange's work was used in geography.  The article contains details on the choice of the constants. 

\section{In guise of a conclusion}

Let us conclude this chapter by recalling that Lagrange, besides being a mathematician, was also an astronomer,\footnote{For instance, Lagrange's work on the two body problem\index{two body problem} (applied to the Eath-moon, to the Earth-sun, and to other pairs), in which he introduced the equilibrium points known now under the name Lagrange points, is a major contribution to celestial mechanics.} and that he was aware of the fact that astronomical observations, together with  spherical trigonometry, were at the foundations of geodesy. Let me quote from  \S 9 of his memoir \emph{Solution de quelques probl\`emes relatifs aux triangles sph\'eriques, avec une analyse compl\`ete de ces triangles} (Solution of some problems relative to spherical triangles, with a complete analysis of these triangles) \cite{Lagrange-Solution},  mentioning geographical applications of his trigonometric formula. He writes: 
%``[\ldots] la formule que nous venons de donner sera \'egalement utile pour mesurer les surfaces sph\'eriques termin\'ees par des arcs de grands cercles. Ainsi elle peut être
%employ\'ee avec beaucoup d'avantage pour d\'eterminer l'\'etendue d'un pays, lorsqu'on connaît les latitudes et les diff\'erences de longitude de plusieurs points plac\'es \`a la circonf\'erence; car, en liant ces points par des arcs de grands cercles, on aura un polygone sph\'erique, dont on trouvera facilement l'aire en le d\'ecomposant en quadrilat\`eres form\'es par les cercles de latitude et par les arcs de l'\'equateur intercept\'es entre ces cercles."
``[\ldots] the formula we just gave will also be useful for measuring spherical surfaces terminated by arcs of great circles. Thus it can be
employed with great advantage to determine the extent of a country, when we know the latitudes and differences in longitude of several points placed at the circumference; for, by linking these points by arcs of great circles, we shall have a spherical polygon, whose area we shall easily find by decomposing it into quadrilaterals formed by the circles of latitude and by the arcs of the equator intercepted between these circles."

\end{document}